\documentclass[12pt,notitlepage,twoside]{article}
\usepackage{latexsym} 
\usepackage{amsmath} 
\usepackage{amsfonts} 
\usepackage{amssymb} 
\renewcommand{\a }{\alpha } 
 
\renewcommand{\d}{\delta }

\newcommand{\D }{\Delta }
 
\newcommand{\e }{\varepsilon } 
\newcommand{\g }{\gamma} 
 
\renewcommand{\l }{\lambda }

\newcommand{\n }{\nabla } 
\newcommand{\var }{\varphi }

\renewcommand{\o }{\omega } 
\renewcommand{\O }{\Omega }

\newcommand{\wk}{\rightharpoonup} 
\newcommand{\be}{\begin{equation}} 
\newcommand{\ee}{\end{equation}} 
\newenvironment{pf}{\noindent{\bf Proof.}\enspace}{
\hfill$\Box$\medskip} 
\newenvironment{pfn}[1]{\noindent{\bf Proof of {#1}\enspace}}{
\hfill$\Box$\medskip} 
\newcommand{\R}{\mathbb{R}} 
 
\newtheorem{thm}{Theorem}[section] 
\newtheorem{pro}[thm]{Proposition}
\newtheorem{lem}[thm]{Lemma}
\newtheorem{rem}[thm]{Remark}
\newtheorem{cor}[thm]{Corollary}
 
\numberwithin{equation}{section}

\author{ Khalil El Mehdi$^a$ \& Massimo Grossi$^b$\footnote{Supported by M.U.R.S.T., project
: ``Variational methods and nonlinear differential equations''.}\\
{\footnotesize
a :  Facult\'e des Sciences et Techniques, Universit\'e de Nouakchott, BP 5026, Nouakchott,}\\
{\footnotesize
 Mauritania. E-mail : \texttt{khalil@univ-nkc.mr}}\\
{\footnotesize
 b :  Dipartimento di Matematica, Universit{\'a} di
 Roma "La Sapienza", P. le A. Moro 2,}\\
{\footnotesize
 00185 Roma, Italy. E-mail: \texttt{grossi@mat.uniroma1.it}}
}

\title { \Large \textbf{Asymptotic Estimates and Qualitative
    Properties \\
of an Elliptic Problem in Dimension Two }} 
 
\begin{document} 
 
\date{ } 
 
\maketitle

{\footnotesize 

\noindent 
{\bf Abstract.-}
 In this paper we study a semilinear elliptic problem on a bounded
 domain in $\R ^2$ with large exponent in the nonlinear term. We
 consider positive solutions obtained by minimizing suitable
 functionals. We prove some asymptotic estimates which enable us to
 associate a "limit problem" to the initial one. Using these estimates we prove some qualitative
properties of the solution, namely characterization of  level sets and  nondegeneracy.
 
\noindent\footnotesize {{\bf Classification AMS:}\quad 35J60}

}

\section{Introduction and Main Results }
In this paper we consider the following elliptic problem 
 $$
( P_{\l,p}) \quad \left\{
\begin{array}{cccc}
 -\Delta u+\l u =u^{p} & \mbox{ in }\, \O \\ 
         u>0             & \mbox{  in  }\,  \Omega\\
         u=0             & \mbox{ on }\,  \partial  \Omega 
\end{array}
\right.
$$ 
\noindent
where $\Omega$ is a bounded domain in $\R^2 $,  $ \l \geq 0 $
and $ p $ is a large positive parameter.\\
We will focus on the solutions to $ ( P_{\l,p} )$ obtained by the following
variational method.\\
We define on $ H^1_0 (\O) \diagdown \{0\} $ the  $C^2$-functional
$$
J_\l (u) = \frac {\int _\O |\nabla u|^2 + \l \int _\O u^2 } 
{(\int _\O |u|^{p+1})^{2/(p+1)}}
$$
and we consider the following minimizing problem
\begin{equation}
\label{e:11} 
c_{\l,p}^2:=\inf _{u\in H^1_0(\O) \diagdown \{0\}} J_\l (u). 
\end{equation}
\noindent
A standard variational argument shows that $ c_{\l,p}^2 $ can be achieved
by a positive function. Then after a multiplicative constant we find a
positive function $  u_{\l,p} $ which solves $ ( P_{\l,p} ) $ and
satisfies
\begin{equation}
\label{e:111} 
c_{\l,p}^2= \frac {\int _\O |\nabla  u_{\l,p}|^2 + \l \int _\O  u_{\l,p}^2} 
{(\int _\O | u_{\l,p}|^{p+1})^{2/(p+1)}}.
\end{equation}

In the remainder of this paper we denote by $  u_{\l,p} $ a least energy
solution of $ ( P_{\l,p}) $ obtained in this way.\\

The aim of this paper is to study qualitative properties of the solution $u_{\l,p}$ for $\l\ge0$ and $p$
large. An essential tool in the proof of these results is to have information on the asymptotic
behavior of $ u_{\l,p} $ as $p$ becomes large. 
The asymptotic behavior of the solutions of $P_{\l,p}$ was initially studied by Ren and Wei when
$\l=0$. More precisely, in \cite {RW1} and \cite {RW2} the authors proved 
the following result:
\begin{thm}(\cite {RW1}, \cite {RW2})\label{t:11}
Let $ \O $ be a smooth bounded domain and $ \l = 0 $ in $ ( P_{\l,p} )
$. Let us denote by $  u_{0,p} $ a least energy solution of $ J_0
(u)$. Then, for any sequence $u_{p_n}$ of $ u_{0,p}$ with $p_n \to +\infty
$, there exists a subsequence of $u_{p_n}$, still denoted by $u_{p_n}$,
such that\\
\textbf{i.} $ u_{p_n}^{p_n} (\int_\O u_{p_n}^{p_n})^{-1} \to \d _{x_0} $ in the
sense of distribution, where $\d_{x_0} $ is the Dirac function at point
$x_0$.\\
\textbf{ii.} $ u_{p_n}^{p_n} (\int_\O u_{p_n}^{p_n})^{-1} \to G(x, x_0
) $ in $W^{1,q}(\O ) $ weakly for any $ 1< q < 2 $, where $G $ is the
Green's function of $-\D $ with Dirichlet boundary
condition. Furthermore, for any compact set $K \subset \bar{\O}
\diagdown \{x_0\} $, we have $v_{p_n} \to G(., x_0 ) $ in $C^{2, \a
}(K) $.\\
\textbf{iii.} $x_0 $ is a critical point of the Robin function $R $ defined by 
$R(x) = g(x,x) $, where 
$$
 g(x,y) = G(x,y) + \frac{1}{2\pi}Log|x-y|
$$
is the regular part of the Green's function.  
\end{thm}
\noindent
Moreover in \cite{RW1} it was also showed that
\begin{eqnarray}\label{e :12}
0 < C_1 \leq || u_{0,p} ||_{L^\infty (\O )} \leq C_2
\end{eqnarray}
for some constants $C_1 $, $C_2 $ and for $p $ large enough. From
these results we can see that when $p $ gets large, the least energy  solution
$ u_{0,p} $ looks like a single spike. \\

One of the results of this paper is to obtain asymptotic estimates for the least
energy solution $ u_{\l,p} $, but of different type of to the
corresponding one due to Ren and Wei. To describe our results we need to
introduce the following problem
\begin{eqnarray}\label{e:13}
\left\{
\begin{array}{cccc}
-\D u = e^u & \mbox {in } \R^2 \\
 \int _{\R^2} e^u<+\infty&
\end{array}
\right.
\end{eqnarray}
In \cite{CL} it was proved that any solution of \eqref{e:13} is given
by
\begin{eqnarray}\label{e:14}
U_{\mu , y}(x) = Log\left(\frac {8 \mu ^2}{(1+\mu ^2|x-y|^2)^2}\right)
\end{eqnarray}
with $ \mu \in \R $ and $ y \in \R^2 $. \\
Now we can claim the following
\begin{thm}\label{t:12}
Let $ \O $ be a smooth bounded domain of $\R^2 $, $\l\ge0$ and let $ u_{\l,p} $ be a
least energy solution of $( P_{\l,p} ) $. Then, we have \\
\textbf{i.} $|| u_{\l,p}||_\infty ^{p-1} \to +\infty $ as $ p \to +\infty$.
\\
\textbf{ii.}  If $ \var_{\l,p} $ is the function defined for $x \in \O_{\l,p} := 
|| u_{\l,p}||_\infty ^{(p-1)/2}(\O -x_{\l,p})$
$$
\var_{\l,p} (x) = (p-1)Log\Big(\frac{u_{\l,p}}{|| u_{\l,p}||_\infty}\big(x_{\l,p} +\frac{x}
{\sqrt{p-1} || u_{\l,p}||_\infty ^{(p-1)/2 }}\big)\Big)
$$
where $ x_{\l,p} \in \O $ is such that $ || u_{\l,p}||_\infty =  u_{\l,p} (x_{\l,p} ) $, then,
for any sequence $\var _{\l ,p_n} $ of $\var _{\l ,p} $ with $ p_n \to \infty $,
there exists a subsequence of $\var _{\l ,p_n} $, still denoted by $\var
_{\l ,p_n} $, such that  $ \var _{\l ,p_n} \to U_{\bar{\mu}, 0} $ in $
C^2_{loc}(\R^2) $, where $ \bar{\mu}^2 = 1/8 $ and $
U_{\bar{\mu},0}$ is given by \eqref{e:14}.
\end{thm}
Since $|| u_{\l,p}||_\infty ^{(p-1)/2} \to \infty $ and $ x_{\l ,p} \to x_0 \in \O
$ (see Corollary \ref{c:25} below), we have that $ \O _{\l ,p} \to \R^2 $ as $ p \to
\infty $. From this, we say that \eqref{e:13} is the "limit problem "
of $( P_{\l,p} ) $ as $ p \to \infty $.

A similar phenomenon (existence of
a "limit problem ") occurs in several
situations in higher dimensions. A typical example is the following
problem
\begin{eqnarray}\label{e:15}
\left\{
\begin{array}{cccc}
-\D u+\l u = n(n-2)u^{\frac{n+2}{n-2} - \e }& \mbox { in } \O\\
    u > 0                              &\mbox { in } \O\\
    u = 0                              &\mbox {on } \O
\end{array}
\right.
\end{eqnarray}
where $\e $ is a small positive parameter and $n\geq 3 $. Here it is
well known that the limit problem associated to \eqref{e:15} is
\begin{eqnarray}\label{e:16}
\left\{
\begin{array}{cccc}
-\D u = n(n-2)u^{\frac{n+2}{n-2} }& \mbox { in } \R^n\\
    u > 0                              &\mbox { in } \R^n
\end{array}
\right.
\end{eqnarray}  
which admits the two parameters family of solutions
$$
\d _{\mu , y} (x) = \frac{\mu ^{(n-2)/2}}{(1+\mu ^2|x-y|^2)^{(n-2)/2}}.
$$
Theorem \ref{t:12} emphasizes some similarities between the problem $
( P_{\l,p} ) $ when $p $ is large and some corresponding problems in higher
dimensions. 

We remark that Theorem \ref{t:12} is the starting
point to obtain similar results as in singularity perturbed problems
involving the critical Sobolev exponent, namely uniqueness or
qualitative properties of solutions. Proof of Theorem \ref{t:12} is given in Section 2.

In Section 3 we give a first application of Theorem
\ref{t:12}; we study the shape of the level sets
of solutions $ u_{\l,p} $ when $ p $ is large enough. Namely we have the
following result.
\begin{thm}\label{t:13}
Let $u_{\l,p}$ be a least energy solution to $P_{\l,p}$ satisfying (\ref{e:11}). Let $\O$ be
convex.

Then there exists $p_0\geq 1$  such that for any $p>p_0$, we have
$$
(x-x_{p})\nabla  u_{\l,p} (x) < 0, \quad \forall x \in \O
\diagdown \{x_{p}\}
$$
where $x_{p} \in \O $ such that $  u_{\l,p} (x_{p} ) = || u_{\l,p} ||_\infty $.\\
In particular, $x_{p} $ is the only critical point and the superlevels
are strictly star shaped with respect to $ x_{p} $ for $p $ large enough.
\end{thm}

If $\O$ is also symmetric, the claim of Theorem \ref{t:13} follows by the well known Gidas-Ni-Nirenberg Theorem.
Notice that Theorem \ref{t:13} was proved by Lin (\cite{L}) if $\lambda=0$ and $p>1$ with different techniques.

In Section 4 we give another application of Theorem \ref{t:12}, proving some uniqueness and
nondegeneracy result to $P_{\l,p}$ for domains which satisfy the assumption of the Gidas-Ni-Nirenberg
Theorem. 

\begin{thm}
\label{t:14}
Let $\O$ be a smooth bounded domain of $\R^2$ which is symmetric with respect to the plane $x_1=0$
and $x_2=0$ and convex with respect to the direction $x_1$ and $x_2$. Let $u_{\l,p}$ be a least
energy solution of $P_{\l,p}$. Then there exist $p_0\ge1$  such that for any $p\ge p_0$ we have that $u_{\l,p}$ is nondegenerate, i.e. the problem
\begin{equation}
\label{21}
\left\{
\begin{array}{cccc}
 -\Delta v+\l v =pu_{\l,p}^{p-1}v & \mbox{ in }\, \O \\ 
         v=0             & \mbox{ on }\,  \partial  \Omega 
\end{array}
\right.
\end{equation}
admits only the trivial solution $v\equiv0$.
\end{thm}

Similar ideas used in the proof of Theorem \ref{t:14} could help to obtain uniqueness result for the least energy solution to $(P_{\l,p})$. It will be done in a forthcoming paper.
\section{Proof of Theorem \ref{t:12}}
In this section we give the proof of Theorem \ref{t:12}. 
Here we suppose that $\l > 0$ is fixed. We begin by proving some auxiliary lemmas.
\begin{lem}\label{l:21}
There exists $c>0$ such that $||u_{\l,p}||_\infty \geq c $, where $u_{\l,p}$ is a solution of
$( P_{\l,p} )$ and $c$ is independent of $p$.
\end{lem}
\begin{pf}
Let $\l _1$ be the first eigenvalue of $-\D $ and $e_1 $ be a
corresponding positive eigenfunction. Then if $u_{\l,p}$ is a solution of
$( P_{\l,p} )$, we have
$$
0= \int _\O (u_{\l,p}\D e_1 - e_1\D u_{\l,p}) = -\l_1 \int _\O u_{\l,p}e_1 + \int _\O e_1
(u_{\l,p}^p -\l u_{\l,p}).
$$
Thus 
$$
\int_\O e_1u_{\l,p}^p = (\l +\l _1)\int _\O e_1u_{\l,p}.
$$
Hence
$$
(\l +\l _1)\int _\O e_1 u_{\l,p} \leq ||u_{\l,p}||_\infty ^{p-1} \int _\O e_1 u_{\l,p}.
$$
Then
$$
||u_{\l,p}||_\infty\ge(\l + \l _1)^{\frac{1}{p-1}}\ge\min\{\l_1,1\}.
$$
Therefore our lemma follows.
\end{pf}\\
\noindent
\begin{lem}\label{l:12}
For $p$ large enough, there exists $c$ such that
$$
c_{\l ,p} \leq c\, p^{-1/2}
$$
where $c_{\l ,p}$ is defined in \eqref{e:11}.
\end{lem}
\begin{pf}
We follow the proof of Lemma 2.2 in \cite {RW1}. Without loss of
generality we can assume $0\in \O $. Let $R > 0$ be such that
$B(0,R)\subset \O$. For $0 < d < R$, we introduce the following Moser
function
$$
m_d(x) =\frac{1}{\sqrt{2\pi}} 
\begin{cases}
(Log (R/d))^{1/2}\quad \quad \quad \quad \quad \quad \mbox {if }\, 0\leq |x| \leq d \\
Log (R/|x|)(Log (R/d))^{-1/2}\quad \mbox {if }\, d\leq|x|\leq R\\
0 \quad \quad \quad \quad \quad \quad \quad \quad \quad \quad \quad\quad \mbox {if }\, |x| \geq R
\end{cases}
$$
Then $m_d \in H^1_0(\O)$ and $||\nabla m_d||_{L^2(\O )} = 1$.\\
Observe that
$$
\int _\O m_d^{p+1}(x)dx = I_1 + I_2
$$
where
$$
I_1 = \left(\frac{1}{\sqrt{2\pi}}(Log(R/d))^{1/2}\right)^{p+1}\pi d^2
$$
and
$$
I_2= \left(\frac{1}{\sqrt{2\pi}}(Log(R/d))^{-1/2}\right)^{p+1}\int
_{d<|x|<R}(Log(R/|x|))^{p+1}dx
$$
Thus
$$
|m_d|_{L^{p+1}(\O )}\geq I_1^{1/(p+1)}
$$
Choosing $d = Re^{-(p+1)/4}$, we find
$$
|m_d|_{L^{p+1}}^2 \geq (p+1)(8\pi e)^{-1} (\pi R^2)^{2/(p+1)}
$$
Hence
$$
\frac{\int |\nabla m_d|^2 + \l \int m_d^2}{|m_d|_{L^{p+1}}^2}
\leq \Big(1+\frac{\l'}{\l_1(B(0,R)}\Big)\frac{\int |\nabla m_d|^2}{|m_d|_{L^{p+1}}^2}
\leq c(R)(p+1)^{-1}R^{-4/(p+1)}
$$
Then
$$
c_{\l ,p} \leq c(R)(p+1)^{-1/2}R^{2/(p+1)}
$$
Therefore our lemma follows.
\end{pf}\\
\noindent
\vskip0.2cm
In addition, since our solution $u_{\l,p}$ satisfies (\ref{e:111}) and
$$
\int _\O |\nabla  u_{\l,p} |^2 +\l \int _\O  u_{\l,p}^2 = \int _\O  u_{\l,p}^{p+1}
$$
we easily derive the following result
\begin{cor}\label{c:23}
For $p$ large enough, there exists $c > 0$ such
that
\begin{equation}
\label{b0}
p\int _\O  u_{\l,p}^{p+1} \leq c \quad \mbox{and} \quad 
p(\int _\O |\nabla u_{\l ,p}|^2 + \int _\O u_{\l ,p}^2 ) \leq c
\end{equation}
\end{cor}
\noindent
Now, we recall the following lemma (see \cite{DLN}, \cite{GNN})
\begin{lem}(\cite {DLN} \cite{GNN})\label{l:24}
Let $u$ be a solution of
$$
\left\{
\begin{array}{ccccc}
-\D u = F(u)& \mbox {in }\O \subset \R^2 \\
    u = 0   & \mbox {on } \partial\O
\end{array}
\right.
$$
where $\O$ is a bounded smooth domain and $F$ is a
$C^1$-function. Then, there exists a neighborhood $\o$ of $\partial\O$
and $C > 0$, both depending only on $\O$, such that 
$$
||u||_{L^{\infty}(\omega)} \leq C ||u||_{L^1(\O)}
$$
\end{lem}
\begin{cor}\label{c:25}
Let us denote by $x_{\l,p}$ the point where $u_{\l,p}$ achieves its maximum,
that is $||u_{\l,p}||_\infty = u_{\l,p}(x_{\l,p}) $. Then  $x_{\l,p} \to x_\l \in \O $.
\end{cor}
\begin{pf}
From Lemma \ref{l:21}, we have $||u_{\l,p}||_\infty \geq c > 0 $ and from
Corollary \ref {c:23}, we derive that
$$
\int _\O u_{\l,p} \to 0.
$$
Using Lemma \ref{l:24} we deduce that the point $x_{\l,p}$ is far away
from the boundary. Thus the claim follows.
\end{pf}\\
\begin{lem}\label{l:26}
There exist a sequence $p_n\rightarrow\infty$ such that
$$
\lim\limits_{n\rightarrow\infty}
||u_{\l,p_n}||_\infty ^{p_n-1} \to +\infty.
$$
\end{lem}
\begin{pf}
We argue by contradiction. Let us suppose that there exists $c>0$ such that for any $p>1$ and
$\l > 0$ we have 

\begin{equation}
\label{b1}
|| u_{\l,p}||_\infty ^{p-1}\leq c .
\end{equation}
Let us consider the following function
\begin{eqnarray}
\bar u_{\l,p}(X) = \frac{1}{|| u_{\l,p}||_\infty } u_{\l,p} (x_{\l ,p} +
\frac{X}{|| u_{\l,p}||_\infty ^{(p-1)/2}}) \quad \mbox{for } X\in \O_{\l,p}
\end{eqnarray}
where $\O_{\l,p} = || u_{\l,p}||_\infty ^{(p-1)/2} (\O - x_{\l,p})$ and $x_{\l,p} \in \O $
such that $ u_{\l,p}(x_{\l,p}) = || u_{\l,p}||_\infty $.\\
\noindent
It is easy to see that $\bar{u}_{\l ,p}$ satisfies
\begin{eqnarray}
\left\{
\begin{array}{cccc}
-\D \bar u_{\l,p} = \bar u_{\l,p}^p - \frac{\l }{|| u_{\l,p}||_\infty ^{p-1}}
 \bar u_{\l,p} & \mbox {in } \O_{\l,p} \\
 \bar u_{\l,p}(0) =1,\quad 0\leq \bar u_{\l,p}\leq 1 & \mbox {in }\O_{\l,p}\\
  \bar u_{\l,p} =0 & \mbox {on } \partial \O_{\l,p}
\end{array}
\right.
\end{eqnarray}
Thus by the standard regularity theory we deduce that there exists a sequence
$p_n\rightarrow\infty$ such that $\bar u_{\l ,p_n}\rightarrow\bar{u}_\l$ in $C^1_{loc}(\R^2)$. Moreover by
(\ref{b1}), up to a subsequence of $p_n$,
$\O_{\l ,p_n}\to D:=\g _\l (\O - x_\l ) $ as $n\to +\infty $ with
$\g _\l=\lim\limits_{n\to\infty}||u_{\l ,p_n}||_\infty^{p_n-1}$. Let us point out that, by
$0\le\bar{u}_{\l ,p}\le1$, we derive that $\bar{u}_{\l
  ,p_n}^{p_n}\wk\psi _\l\ge0$ weakly in $L^q(D)$ for any
$q>1$.
 Finally $\bar{u}_\l$ satisfies
\begin{eqnarray*}
\left\{
\begin{array}{cccc}
-\D \bar{u}_\l = \psi _\l - \bar{\l}\bar{u}_\l & \mbox {in } D \\
 \bar{u}_\l (0) =1,\quad 0\leq \bar{u}_\l \leq 1 & \mbox {in }D\\
  \bar{u}_\l =0 & \mbox {on } \partial D
\end{array}
\right.
\end{eqnarray*}
where $\bar{\l} = \lim\limits_{n\to +\infty}\l ||u_{\l ,p_n}||_\infty
^{1-{p_n}}=\frac{\l}{\g _\l} $.\\
Thus
$$
\int _D |\nabla \bar{u}_\l|^2 + \bar{\l} \int _D \bar{u}_\l ^2 = \int _D
\psi _\l \bar{u}_\l
$$
Observe that, by Lebesgue's Theorem and the definition of $\bar{u}_\l$,
we have
\begin{eqnarray*}
\int _D \psi _\l \bar{u}_\l = \lim _{n\to \infty }\int _{\O_{\l
    ,p_n}}\bar{u}_{\l ,p_n}^{p_n}. 
\bar{u}_{\l ,p_n} = \lim _{n\to \infty}\int _\O ||u_{\l ,p_n}||_\infty ^{-2}
u_{\l ,p_n}^{{p_n}+1} 
\end{eqnarray*}
From Corollary \ref{c:23} and Lemma \ref{l:21}, we derive
$$
\int _D \psi _\l \bar{u}_\l=0
$$
Hence
$$
\int _D |\nabla \bar{u}_\l|^2 + \bar{\l} \int _D \bar{u}_\l^2 = 0
$$
Therefore $\bar{u}_\l \equiv 0 $ which gives a contradiction with
$\bar{u}_\l(0)=1 $ and our lemma follows.
\end{pf}\\
\begin{lem}\label{l:27}
For any $ x\in B(x_{\l ,p}, \frac{R}{\sqrt{p-1}|| u_{\l,p}||_\infty ^{(p-1)/2}})
$ we have, for $p$ large enough,
$$
 u_{\l,p}(x) \geq \g _R > 0
$$
where $R$ is an arbitrary positive number and $\g _R $ is a constant
only depending on $R$.
\end{lem}
\begin{pf}
For $X\in \O_{\l,p}$, we set 
$$
W_{\l,p}(X) = u_{\l,p}(x_{p} + \frac{X}{\sqrt{p-1}|| u_{\l,p}||_\infty ^{(p-1)/2}}).
$$
Thus $W_{\l,p}$ satisfies
$$
-\D W_{\l,p} (X) = c_{\l,p}(X) W_{\l,p}(X)
$$
where 
$$
c_{\l,p}(X) = \frac{u_{\l,p}^{p-1}}{(p-1)|| u_{\l,p}||_\infty ^{p-1}}(x_{p} +
\frac{X}{\sqrt{p-1}|| u_{\l,p}||_\infty ^{p-1}}) -
\frac{\l}{(p-1)|| u_{\l,p}||_\infty ^{p-1}}
$$
Observe that, from Lemma \ref{l:26}, we deduce
$$
|c_{\l ,p}(X)|\leq C \quad \quad \mbox{for $p$ large enough}.
$$
From the standard Harnack inequality \cite{GT}, we get
$$
|| u_{\l,p}||_\infty = \sup _{B(0,R)} W_{\l ,p}(X) \leq c_R \inf
_{B(0,R)}W_{\l ,p}(X)
\quad \mbox{for $p$ large enough.}
$$
Thus
$$
\inf _{B(0,R)}W_{\l ,p} \geq \frac{|| u_{\l,p}||_\infty }{c_R}
$$
From Lemma \ref{l:21}, we deduce
$$
\inf _{B(0,R)}W_{\l ,p} \geq \g _R
$$
and therefore our lemma follows. 
\end{pf}\\
\begin{lem}\label{l:28} 
Let us consider the function
\begin{eqnarray}\label{e:24} 
F_{\l,p}(X)= \frac{1}{|| u_{\l,p}||_\infty ^{p-1}}\frac{|\nabla  u_{\l,p}(x_{\l,p} + \frac
  {X}{\sqrt{p-1}|| u_{\l,p}||_\infty ^{(p-1)/2}})|^2}{ u_{\l,p}^2(x_{\l,p} + \frac
  {X}{\sqrt{p-1}|| u_{\l,p}||_\infty ^{(p-1)/2}})}, \quad \mbox{for } X\in
  \O _{\l,p}
\end{eqnarray}
Then, for any $R > 0 $, we have, for $p$ large enough
$$
||F_{\l,p}||_{L^\infty (B(0,R))} \leq C_R
$$
where $ C_R $ is a constant only depending on $ R$.
\end{lem}
\begin{pf}
According to Lemma \ref{l:27}, it is enough to prove 
$$
 \frac{1}{|| u_{\l,p}||_\infty ^{p-1}}\Big|\nabla  u_{\l,p}\Big(x_{p} + \frac
  {X}{\sqrt{p-1}|| u_{\l,p}||_\infty ^{(p-1)/2}}\Big)\Big|^2 \leq c
$$
For $ X\in B(0,2R) $, let 
$$
f_{\l,p}^i(X) = \frac{1}{|| u_{\l,p}||_\infty ^{(p-1)/2}}\frac{\partial
   u_{\l,p}}{\partial x_i}\Big(x_{p} +\frac
  {X}{\sqrt{p-1}|| u_{\l,p}||_\infty ^{(p-1)/2}}\Big), \quad \mbox {for } i=1, 2
$$
It is sufficient to prove
$$
|f_{\l,p}^i|_{L^\infty (B(0,R))} \leq c, \quad \mbox{for }i=1, 2 \mbox{ and
  c is independent of}\ p.
$$
We point out that
$$
-\D f_{\l,p}^i = c_{\l,p} (X) f_{\l,p}^i
$$
with 
$$
c_{\l,p}(X)= \frac{p}{p-1}\frac{u_{\l,p}^{p-1}}{|| u_{\l,p}||_\infty ^{p-1}} (x_{p} + \frac
  {X}{\sqrt{p-1}|| u_{\l,p}||_\infty ^{(p-1)/2}}) -
  \frac{\l}{(p-1)|| u_{\l,p}||_\infty ^{p-1}}
$$
From Lemma \ref{l:26} we have
$$
|c_{\l,p}(X)|\leq C \quad \quad \mbox{for $p$ large enough.}
$$
Hence, by the standard weak Harnack inequalities (Theorem 8.17 of
\cite{GT}), we have
$$
||f_{\l,p}^i||_{L^\infty (B(0,R))} \leq c ||f_{\l,p}^i||_{L^2(B(0,R))}.
$$
Observe that
\begin{align*}
||f_{\l,p}^i||^2_{L^2(B(0,2R))}&= \int _{B(0,2R)}\frac{1}{|| u_{\l,p}||_\infty
  ^{p-1}}
|\frac{\partial  u_{\l,p}}{\partial x_i}(x_{p} + \frac
  {X}{\sqrt{p-1}|| u_{\l,p}||_\infty ^{(p-1)/2}})|^2 dX\\
&=(p-1)\int _{B(x_{\l,p}, \frac{2R}{\sqrt{p-1}|| u_{\l,p}||_\infty ^{(p-1)/2}})}
\Big|\frac{\partial  u_{\l,p}}{\partial x_i}(x)\Big|^2 dx\\
&\leq  (p-1)\int _\O |\nabla  u_{\l,p}|^2.
\end{align*}
From Corollary \ref{c:23}, we derive
$$
||f_{\l,p}^i||_{L^2(B(0,2R))} \leq c.
$$
Therefore our lemma follows
\end{pf}\\ 
Next we will prove Theorem \ref{t:12}.\\
\begin{pfn}{\bf Theorem \ref{t:12}}  
According to Lemma \ref{l:26}, it only remains to prove part ii. of
the Theorem. To do this we introduce the following function
\begin{eqnarray}\label{e:25}
v_{\l,p}(X)= \var _{\l,p}(X) - Z_{\l,p}(X), \quad \mbox{for } X\in \O_{\l,p}
\end{eqnarray}
where $Z_{\l,p}$ satisfies
\begin{eqnarray}\label{e:26}
\left\{
\begin{array}{ccccc}
-\D Z_{\l,p}+\frac{\l}{|| u_{\l,p}||_\infty ^{p-1}}= F_{\l,p}&\mbox{in }&B(0,R)\\
    Z_{\l,p}=0                            &\mbox{on }&\partial B(0,R)
\end{array}
\right.
\end{eqnarray}
where $\var_{\l,p}$ is defined in Theorem \ref{t:12} and $F_{\l,p}$ is
defined by \eqref{e:24}. By the maximum principle we have that $Z_{\l,p}\ge0$
From Lemmas \ref{l:26} and \ref{l:28} and the standard regularity
theory, we derive that for any $R>0 $
\begin{eqnarray}\label{e:27}
||Z_{\l,p}||_{C^1(B(0,R))} \leq C_R
\end{eqnarray}
where $C_R $ only depends on $R$.\\
Thus setting
\begin{eqnarray}\label{e:28}
V_{\l,p}(X) = e^{Z_{\l,p}(X)}
\end{eqnarray}
we have that, for any $q\ge1$,
$$
\forall R > 0, \, \, \exists C_R > 0 \quad\hbox{such that }||V_{\l,p}||_{L^q(B(0,R))} \leq
C_R.
$$
By direct computation it is not difficult to see that $v_{\l,p}$ satisfies
\begin{eqnarray}\label{M}
\left\{
\begin{array}{ccccc}
-\D v_{\l,p} = V_{\l,p}(x) e^{v_{\l,p}}& \mbox{in }B(0,R) \\
    v_{\l,p} \leq  0      &\mbox{in }B(0,R)
\end{array}
\right.
\end{eqnarray}     
We claim that\\
{\bf a.} $ V_{\l,p} \geq 0 $ in $ B(0,R) $ \\
{\bf b.} $ ||V_{\l,p}||_{L^q(B(0,R))} \leq C_R  \quad \forall q\ge1$ \\
{\bf c.} $ \int _{B(0,R)}e^{qv_{\l,p}} \leq C'_R \quad \forall q\ge1$\\
Note that {\bf a.} and {\bf b.} follow by the definition of $ V_{\l,p}(x)$. Concerning {\bf c.} we have
that
\begin{eqnarray}\label{z}
& &\int _{B(0,R)}e^{qv_{\l,p}}\le\int _{B(0,R)}e^{q\var _{\l,p}}=
\int _{B(0,R)}\Big[\frac{ u_{\l,p}^{p-1}}{|| u_{\l,p}^{p-1}||_{\infty}^{p-1}}
\big(x_{\l,p} + \frac{X}{(p-1)^{1/2}
|| u_{\l,p}||_\infty ^{(p-1)/2 }}\Big]^q\nonumber\\
& &\le\int _{\Omega_{\l,p}}\frac{ u_{\l,p}^{p-1}}{|| u_{\l,p}^{p-1}||_{\infty}^{p-1}}
\Big(x_{\l,p} + \frac{X}{(p-1)^{1/2}
|| u_{\l,p}||_\infty ^{(p-1)/2 }}\Big)=(p-1)\int _{\Omega} u_{\l,p}^{p-1}\le C
\end{eqnarray}

Thus, we are in the setting of Theorem 3 of Brezis-Merle \cite{BM} and
we then have the following alternative :\\
either\\
(i) $v_{\l,p} $ is bounded in $L^\infty (B(0,R))$\\
or\\
(ii) $ v_{\l,p} \to - \infty $ uniformly in $B(0,R) $\\
or\\
(iii) $ v_{\l,p} \to -\infty $ uniformly in $B(0,R) \diagdown S $, where
$S$ is the blow-up set of $v_{\l,p}$, i.e. 
\\
$S= \{x\in B(0,R)$ such that there exists a sequence $y_{\l,p} \in B(0,R)$ with $y_{\l,p} \to x
$ and $v_{\l,p}(y_{\l,p}) \to +\infty \}.$
\\
Since $v_{\l,p} \leq 0 $, we derive $S=\emptyset $ and so (iii) does not
occur. Let us also prove that (ii) cannot happen. From \eqref{e:27} it
is sufficient to prove 
$$
\min _{B(0,R)}\var _{\l,p}(X) \geq -C_R
$$
Let us introduce the following function
\begin{eqnarray}\label{e:29}
\psi _{\l,p} (X)= \frac{u_{\l,p}^{p-1}}{|| u_{\l,p}||_\infty ^{p-1}}\Big(x_{\l,p} +
\frac{X}{\sqrt{p-1}|| u_{\l,p}||_\infty ^{(p-1)/2}}\Big), \quad \mbox{for } X\in
\O _{\l,p}
\end{eqnarray}
It is easy to see that $\psi _{\l,p} $ satisfies
$$
-\D \psi _{\l,p} = \psi _{\l,p}^2 - \frac{\l}{|| u_{\l,p}||_\infty ^{p-1}}\psi _{\l,p} - 
\frac{p-2}{p-1}\, \, \frac{|\nabla \psi _{\l,p} |^2}{\psi _{\l,p}}
$$
Hence
$$
-\D \psi _{\l,p}(X) \leq a_{\l,p}(X)\psi _{\l,p}(X)
$$
with $ a_{\l,p}(X)=\psi _{\l,p} (X)-\frac{\lambda}{|| u_{\l,p}||_\infty ^{p-1}} \in (-1,1] $.\\
By standard weak Harnack inequality (see Theorem 8.17 of \cite {GT}), we
derive
$$
1= \sup _{B(0,\frac{R}{2})}\psi _{\l,p}(X) \leq C_R\left(\int _{B(0,R)}\psi _{\l,p}^2\right)^{1/2}
$$
Thus
$$
\int _{B(0,R)}\psi _{\l,p}^2 \geq C_R^{-2}
$$
Hence
$$
\int _{B(0,R)}e^{2\var _{\l,p}} \geq C_R^{-2}
$$
So (ii) also cannot occur. 

Therefore $v_{\l,p} $ is bounded in $L^\infty
(B(0,R)) $. Then $\var _{\l,p}$ is also bounded in $L^\infty(B(0,R))$ since
$$
-C_R \leq \var _{\l,p} = v_{\l,p} + Z_{\l,p} \leq 0 \quad \mbox{in } B(0,R)
$$
Using the standard regularity theory, since $||Z_{\l,p}||_{L^\infty (B(0,R))}\leq C $ and $v_{\l,p} $ is
bounded we derive from (\ref{e:26}) and (\ref{M}) that $Z_{\l,p}$ and $v_{\l,p}$ are both bounded in
$C^1(B(0,R))$. Thus $||\var _{\l,p} ||_{C^1(B(0,R))} \leq C $.\\
We note that $\var_{\l,p} $ satisfies
\begin{eqnarray}\label{e:210}
-\D \var _{\l,p} = -\frac{\l}{|| u_{\l,p}||_\infty ^{p-1}} + \frac{1}{p-1}|\nabla
 \var _{\l,p} |^2 + e^{\var _{\l,p}}
\end{eqnarray}
Again by the standard regularity theory we get $||\var_{\l,p}||_{C^2_{loc}(\R^2)}\le C$.
Then, for any sequence $p_n\to\infty$ there
exists a subsequence (denoted again by $p_n$) such that
${\var}_{\l,p_n}\to\var$ in $C^1_{loc}(\R^2)$.

Let us show that $e^{\var_{\l,p_n}}\to e^{\var}$ in $ L^1_{loc}(\R^2)$. Since 
$ \var_{\l,p_n}
\rightharpoonup \var $ in $ H^1_{loc}(\R^2) $ we have
\begin{eqnarray*} 
\int _{B(0,R)}|e^{\var_{\l,p_n}} - e^{\var}| &=&
\int _{B(0,R)}|\int _0^1 e^{t\var_{\l,p_n} + (1-t)\var}dt||\var_{\l,p_n}
-\var |\\
&\leq & \int _ {B(0,R)}e^{\var_{\l,p_n} + \var } |\var_{\l,p_n} -
\var|\\
&\leq & \left(\int _{B(0,R)}e^{2(\var_{\l,p_n} + \var)}\right)^{1/2}
\left(\int _{B(0,R)}|\var_{\l,p_n} - \var |^2\right)^{1/2}
\end{eqnarray*}
and the claim follows since $ |\var_{\l,p_n} | \leq C $ in $ B(0,R) $. \\
We also note that, from Corollary \ref{c:23}
\begin{eqnarray*}
\lim _{n\to \infty}\int _{B(0,R)} e^{\var_{\l,p_n}} &=& \lim _{n\to +\infty} \int _{B(0,R)}
\frac{1}{||u_{\l,p_n}||_\infty ^{p_n-1}}u_{\l,p_n}^{p_n-1}(x_{p_n} +
\frac{X}{\sqrt{{p_n}-1}||u_{\l,p_n}||_\infty ^{({p_n}-1)/2}})dx\\
&=& \lim _{n\to \infty} ({p_n}-1)\int
_{B(x_{\l,p_n},\frac{R}{\sqrt{{p_n}-1}||u_{\l,p_n}||_\infty ^{({p_n}-1)/2}})} u_{p_n}^{{p_n}-1}(x)dx\\
&\leq & \lim _{n\to \infty }({p_n}-1) \int _\O u_{\l,p_n}^{{p_n}-1}(x)dx\\
&\leq & C
\end{eqnarray*}
where $C$ does not depend on $R$.\\
Then from Fatou's Lemma, we derive
$$
\int _{\R^2}e^{\var }\le\lim\limits_{n\to\infty}\int _{B(O,R)}e^{\var_{\l,p_n}}\le C
$$
Passing to the limit in \eqref{e:210} and using Lemma \ref{l:26}, we
deduce that $\var $ satisfies
\begin{eqnarray*}
\left\{
\begin{array}{cccc}
-\D \var =e^{\var}&\mbox{in } \, \, \R^2\\
 \var (0)=0,\, \, \var \leq 0&\mbox{in }\, \, \R^2\\
\int _{\R^2}e^{\var }dx<\infty&
\end{array}
\right.
\end{eqnarray*}
According to Chen-Li \cite{CL}, we derive $
\var = U_{\bar{\mu},0}
$ ,
where $U_{\bar{\mu},0} $ is defined in \eqref{e:14}.\\
Then our proposition follows.
\end{pfn}
\section{Proof of Theorem \ref{t:13}}
\mbox{}
Let us start by recalling the following result which is a particular
case of a general theorem due to Grossi-Molle \cite{GM}.
\begin{thm}\label{t:31}
Let $\O$ be a smooth domain in $\R^n$, with $n\geq 1$, and $f\in
C^1(\O,\R^+)$. Suppose that $u\in C^3(\O)\cap C^1(\bar{\O})$ satisfies
\begin{eqnarray*}
\left\{
\begin{array}{ccccc}
-\D u+\l u=u^p& \mbox{in}&\O \\
    u >0         &\mbox{in}&\O  \\
    u =0         & \mbox{on}&\partial \O
\end{array}
\right.
\end{eqnarray*}
for $p > 1$ and $\l \in \R$. Let $x_0$ be a maximum point of $u$ and
assume that $\O$ is convex. If there exists an open set $W \subset \O$
containing $x_0$ such that\\
{\bf i.} $(x-x_0) \nabla u(x) < 0, \quad \forall x \in W\diagdown
\{x_0\}$\\
{\bf ii.} $(x-x_0) \nabla u(x) + \frac{2}{p-1} u(x) < 0, \quad \forall x \in
\partial W $\\
{\bf iii.} $\l _1(-\D - pu^{p-1}) > 0 $ in $H^1_0(\O \diagdown W) $(
$\l _1$ is the first eigenvalue of $-\D - pu^{p-1}$)\\
then
$$
(x-x_0) \nabla u(x) < 0 \quad x \in \O \diagdown \{x_0\}
$$
In particular, $x_0$ is the only critical point for $u$ in $\O$ and
the superlevel sets are strictly star shaped with respect to $x_0$.
\end{thm}
\noindent  
For sake of completeness, we recall the proof of Theorem \ref{t:31}.\\
\begin{pfn}{\bf Theorem \ref{t:31}}
 Arguing by contradiction let us suppose that there exists $\bar{x}
 \in \O \diagdown \{x_0\} $ such that $ (\bar{x}-x_0) \nabla u(\bar{x}) \geq
 0 $. By assumption i. $\bar{x} \notin W $. Let us consider
$$
w(x) = (x-x_0)\nabla u(x) +(2/(p-1))u(x)
$$
It turns out that $w(\bar{x}) \geq 0$.\\
Now let us call $D$ the connected component of the set $ \{ x\in \O |
w(x) > 0 \} $ containing $\bar{x}$. By assumption ii. $w < 0 $ on $
\partial W $ and so $ W \cap \partial D = \emptyset $. Moreover if
$z\in \partial \O $,  we have
$$
w(z) = (z-x_0)\nu (z) \frac{\partial u}{\partial \nu}(z)
$$
Since $ \frac{\partial u}{\partial \nu}(z) < 0 $ and using the convexity
of $\O$ we deduce that $ w \leq 0 $ on $\partial \O $. Thus $ w \in
H^1_0(D)$.\\
Now, it is easy to see that $w$ satisfies the following equation
\begin{eqnarray}\label{e :31}
\left\{
\begin{array}{ccc}
-\D w -(pu^ {p-1}+\lambda)w =-2\l u \leq 0 \quad \mbox{in } D\\
w\in H^1_0(D)     & &
\end{array}
\right.
\end{eqnarray}
Since $ \l _1 (-\D -pu^{p-1}I) > 0 $ in $ H^1_0(\O \diagdown W) $ and
$ D \subset \O \diagdown W $,  we get $\l _1(-\D-pu^{p-1}I) > 0 $ in $
H^1_0(D) $. This implies that the maximum principle holds in
$D$. Hence by \eqref{e :31}, we have that $ w \leq 0 $ in $D$ and this
gives a contradiction.
\end{pfn}\\
 In order to prove Theorem \ref{t:13}, we will apply Theorem
 \ref{t:31}.
 Thus we only need to check that the assumptions of Theorem
 \ref{t:31} are true. Let us start by proving the following result.
\begin{pro}\label{p :32}
For any $R > 0 $, we have
$$
(x-x_{\l,p}) \nabla  u_{\l,p} (x) <0, \quad \forall x\in B(x_{\l,p},
\frac{R}{\sqrt{p-1}|| u_{\l,p}||_\infty ^{(p-1)/2}}) \diagdown \{x_{\l,p}\}
$$
for $p$ large enough.
\end{pro}
\begin{pf}
For $x\in B(x_{\l,p},\frac{R}{\sqrt{p-1}|| u_{\l,p}||_\infty
  ^{(p-1)/2}})$ and $ X\in B(0,R) $, we have
$$
(x-x_{\l,p})\cdot\nabla  u_{\l,p}(x) = \frac{1}{\sqrt{p-1}} u_{\l,p}(x)  X\cdot\nabla \var _{\l,p}(X)
$$
where $\var _{\l,p} $ is defined in Theorem \ref{t:12}.\\
Thus it is sufficient to prove that
$$
X\cdot\nabla\var _{\l,p}(X) < 0, \quad \forall X \in B(0,R)\diagdown \{0\}
$$
for $p$ large enough.\\
Arguing by contradiction, let us suppose that there exist $R_0$, a
sequence $ p_n \to +\infty $ and a sequence $\{X_n\}$ in $ B(0,R) $
such that
\begin{eqnarray}\label{e :32}
X_n\cdot\nabla \var _{\l,p_n}(X_n) \geq 0
\end{eqnarray}
From Theorem \ref{t:12}, we know that $ \var _{\l,p_n} \to
U_{\bar{\mu},0} $ in $C^2_{loc}(\R^2)$ where $U_{\bar{\mu},0}$ is
defined \eqref{e:14}. Since $X_n \in B(0,R)$, we can assume that there
exists $X_0\in B(0,2R_0)$ such that $X_n \to X_0$ as $n \to +\infty
$. Thus two cases may occur\\
{\bf Case 1.} $X_0 \not= 0$. Then, in this case, it follows by the above convergence
that
$$
X_n\cdot\nabla \var _{\l,p_n}(X_n) \to X_0\cdot\nabla U_{\bar{\mu},0}(X_0) =
-\frac{4
\bar{\mu}^2 |X_0|^2}{1+ \bar{\mu}^2|X_0|^2} < 0
\quad \mbox{as } n \to +\infty
$$
and this is a contradiction with \eqref{e :32}. Thus this case cannot happen.\\
{\bf Case 2.} $X_0 = 0 $. In this case let us consider the following
function
$$
g_{\l,n}(t)=\var _{\l,p_n}(tX_n), \quad \mbox{for } t \in [0,1]
$$
It yields that $g_{\l,n}$ has a maximum at 0 and another critical point in
[0,1] by \eqref{e :32}(because $ g_{\l,n}'(1) = X_n\cdot\nabla \var _{\l,p_n}(X_n)
\geq 0 $ and $g_{\l,n}'(0)=0 $). Therefore there exists $ \bar{t}_n \in
[0,1] $ such that $ g_{\l,n}''(\bar{t}_n) = 0 $. 
Now let $n\to +\infty $, from the above
convergence and from the assumption $X_0 = 0$, it follows that 0 is a
degenerate critical point for $g_{\l,n}$ and this is not true because $
D^2U_{\bar{\mu},0}(0) = - c Id $, with $c > 0 $. Therefore this
case also cannot happen and our proposition follows.
\end{pf}\\
Now we are able to prove Theorem \ref{t:13}
\vskip0.2cm
\begin{pfn}{\bf Theorem \ref{t:13}}
We will prove that the assumption of Theorem \ref{t:31} are true for
$W=B(x_{\l,p}, \frac{R}{\sqrt{p-1}|| u_{\l,p}||_\infty ^{(p-1)/2}})$.
Proposition \ref{p :32} guarantees that assumption i. holds.\\
Note that
$$
(x-x_{\l,p})\cdot\nabla  u_{\l,p}(x) + \frac{2}{p-1}  u_{\l,p} (x) =
\frac{1}{\sqrt{p-1}}(X\cdot\nabla
\var _{\l,p} (X) + \frac{2}{p-1}) u_{\l,p}(x)
$$
where $\var _{\l,p} $ is defined in Theorem \ref{t:12}.\\
By the convergence of $\var _{\l,p}$ to $U_{\bar{\mu},0}$ and some easy
computations we have that
$$
X\cdot\nabla \var _{\l,p}(X) + \frac{2}{p-1} \to
\frac{-4\bar{\mu}^2|X|^2}{1+\bar{\mu}^2 |X|^2} \quad \mbox{as } p\to
+\infty
$$
uniformly on $ \partial B(0,R)$.\\
Thus
$$
h_{\l,p}(x) := (x-x_{\l,p})\cdot\nabla  u_{\l,p}(x) + \frac{2}{p-1}  u_{\l,p}(x) < 0 \quad \mbox{on }
 \partial B(x_{\l,p}, \frac{R}{\sqrt{p-1}|| u_{\l,p}||_\infty ^{(p-1)/2}})
$$
and then ii. holds.
Finally we have that $h_{\l,p}$ satisfies
\begin{equation}
\label{e :33}
-\D h_{\l,p} - pu^{p-1}h_{\l,p} +\l h_{\l,p} = -2\l  u_{\l,p} < 0
\end{equation}
Since $\nabla  u_{\l,p}(x_{\l,p}) = 0 $, we also have $h_{\l,p}(x_{\l,p}) > 0 $. Thus there
exists a nodal region $C_{\l,p}$ of $h_{\l,p}$ in $B(x_{\l,p},
\frac{R}{\sqrt{p-1}|| u_{\l,p}||_\infty ^{(p-1)/2}})$ where $h_{\l,p} $ is
positive. We derive from (\ref{e :33}) that in $C_{\l,p}$ the first eigenvalue of linearized
operator $L_{\l,p} = -\D -pu^{p-1} +\l $ is negative.
Hence, since $ u_{\l,p}$ is of index 1, the first
eigenvalue of $L_{\l,p}$ in $\O \diagdown C_{\l,p} $ is positive. Thus the first
eigenvalue of $L_{\l,p}$ in $ \O \diagdown B(x_{\l,p},
\frac{R}{\sqrt{p-1}|| u_{\l,p}||_\infty ^{(p-1)/2}})$ is positive. Thus
using Theorem \ref{t:31}, our theorem follows.
\end{pfn}
\vskip0.2cm
\begin{rem}
\label{-3}
It is not difficult to check that Theorem \ref{t:31} holds if $\O$ is star shaped with respect to
$x_0$. Hence, if $\O$ satisfies the assumptions of the Gidas-Ni-Nirenberg Theorem then Theorem
\ref{t:13} holds again.
\end{rem}

\begin{rem}
\label{-4}
It is easy to check that if  $u_{\l,p}$ is a solution to $ P_{\l,p}$, then the first eigenvalue of
the linearized operator $-\Delta +(\l-pu_{\l,p}^{p-1})I_d)$ is negative.

The case where the linearized operator has only nonnegative eigenvalues was considered by various
authors. More precisely, if we consider a solution $u$ of
\begin{equation}
\left\{
\begin{array}{cccc}
 -\Delta u=f(u) & \mbox{ in }\, \O\subset\R^2 \\ 
         u>0             & \mbox{  in  }\,  \Omega\\
         u=0             & \mbox{ on }\,  \partial  \Omega 
\end{array}
\right.
\end{equation}
with $\O$ convex and $\l_1(\Delta-f'(u)I_d)\ge0$ (the so called semistable solution),
Cabr{\'e}-Chanillo (\cite{CC}), Payne (\cite{P}) and Sperb (\cite{S}) showed the uniqueness of the
critical point of $u$.
\end{rem}

\section{A nondegeneracy result}
We start this section recalling the following lemma due to Ren and Wei (\cite{RW1}).
\begin{lem}
\label{-2}\cite{RW1}
For every $t\ge2$ there is $D_t$ such that $||u||_{L^t(\O)}\le D_tt^{\frac{1}{2}}||\nabla
u||_{L^2(\O)}$ for all $u\in H^1_0(\O)$ where $\O$ is a bounded domain in $\R^2$; furthermore
\begin{equation}
\label{-1}
\lim\limits_{t\rightarrow\infty}D_t=(8\pi e)^{-\frac{1}{2}}.
\end{equation}
\end{lem}
From the previous lemma we derive the following estimate, which was showed in \cite{RW1} for $\l=0$.

\begin{lem}
\label{l1} We have, for $p$ large enough and $\l\in[0,\l']$
\begin{equation}
\label{v1} ||u_{\l,p}||_{\infty}\le C
\end{equation} with $C$ depending only on $\l'$.
\end{lem}
\begin{pf}
The proof is the same of the case $\l=0$ (\cite{RW1}), p. 755-756.
Let

\begin{equation}
\label{v2}
\gamma_{\l,p}=\max\limits_{x\in\bar\O}u_{\l,p}(x),\quad{\cal A}=\{x:\frac{\gamma_{\l,p}}{2}
<u_{\l,p}(x)\},\quad\O_t=\{x:t<u_{\l,p}(x)\}.
\end{equation} By Lemma \ref{-2} and Corollary \ref{c:23}

\begin{equation}
\label{v3}
\Big(\int_\O u_{\l,p}^{2p}\Big)^{\frac{1}{2p}}\le C\sqrt p\Big(\int_\O|\nabla u_{\l,p}|^{2}+
\l\int_\O u_{\l,p}^{2}\Big)^{\frac{1}{2}}\le C
\end{equation} for $p$ large and $C$ depending only on $\l''$. Hence

\begin{equation}
\label{v4}
\Big(\frac{\gamma_{\l,p}}{2}\Big)^{2p}|{\cal A}|\le\int_\O u_{\l,p}^{2p}\le C^{2p}.
\end{equation}
 On the other hand
\begin{equation}
\label{v5}
\int_{\O_t}u_{\l,p}^{p}=-\int_{\O_t}\Delta u_{\l,p}+\l\int_{\O_t} u_{\l,p}\ge\int_{\partial\O_t}
|\nabla u_{\l,p}|ds.
\end{equation}
 Using the co-area formula (\cite{F}) and the isoperimetric inequality we have
\begin{equation}
\label{v6} -\frac{d}{dt}|\O_t|\int_{\O_t}u_{\l,p}^{p}\ge\int_{\partial\O_t}\frac{ds}{|\nabla
u_{\l,p}|}
\int_{\partial\O_t}|\nabla u_{\l,p}|ds\ge|\partial\O_t|^2\ge4\pi|\O_t|
\end{equation} From this point we can repeat step by step the proof of (\cite{RW1}), p.756 and we
derive

\begin{equation}
\label{v7}
\gamma_{\l,p}\le C
\end{equation} with $C$ depending only on $\l'$ for $p$ large.
\end{pf}

In the next lemma we study the structure of the solutions of the linearized problem of $P_{\l,p}$
"at infinity". The corresponding result in higher dimensions is well known (\cite{BE},
\cite{R},\cite{AGP}). Here we use some ideas of \cite{AGP} and \cite{D}.
\begin{lem}
\label{0}
Let $v\in L^\infty(\R^2)\cap C^2(\R^2)$ be a solution of the following problem
\begin{equation}
\label{1}
-\Delta v=\frac{1}{(1+\frac{|x|^2}{8})^2}v\quad\hbox{in }\R^2.
\end{equation}
Then
\begin{equation}
\label{2}
v(x)=\sum_{i=1}^2a_i\frac{x_i}{1+\frac{|x|^2}{8}}+b\frac{8-|x|^2}{8+|x|^2}
\end{equation}
\end{lem}
\begin{pf}
We write $v$ as $v=\sum\limits_{k=1}^\infty\psi_k(r)Y_k(\theta)$ where

\begin{equation}
\label{3}
\psi_k(r)=\int_{S^1}v(r,\theta)Y_k(\theta)d\theta,
\end{equation}
and $Y_k(\theta)$ denotes the $k-th$ harmonic spheric satisfying

\begin{equation}
\label{4}
-\Delta_{S^1} Y_k(\theta)=k^2Y_k(\theta)
\end{equation}
Thus (\ref{1}) becomes

\begin{equation}
\label{5}
(-\psi_k''(r)-\frac{1}{r}\psi_k'(r))Y_k(\theta)-\Delta_{S^1} Y_k(\theta)\frac{\psi_k(r)}{r^2}=
\frac{1}{(1+\frac{r^2}{8})^2}Y_k(\theta)\psi_k(r)
\end{equation}
and then
\begin{equation}
\label{6}
-\psi_k''(r)-\frac{1}{r}\psi_k'(r)+k^2\frac{\psi_k(r)}{r^2}=
\frac{1}{(1+\frac{r^2}{8})^2}\psi_k(r)
\end{equation}

Since $v$ is smooth at the origin we deduce that $\psi_k(0)=0$ for $k\ge1$. Moreover since $v\in
L^\infty(\R^2)$ we have that $\psi_k\in L^\infty(\R)$ for any $k\ge0$.

Let us consider the case $k=0$. We have that $\psi_0(r)$ satisfies

\begin{equation}
\label{7}
-\psi_0''(r)-\frac{1}{r}\psi_0'(r)=
\frac{1}{(1+\frac{r^2}{8})^2}\psi_0(r)
\end{equation}
A direct computation shows that $\zeta_0(r)=\frac{8-r^2}{8+r^2}$ is a bounded solution of
(\ref{7}). 
 Let us prove that if $w$ is a second linearly independent solution of (\ref{7}) then
$w$ is not bounded. We write $w(r)=c(r)\zeta_0(r)$. We get from (\ref{7})

\begin{equation}
\label{8}
-(c''\zeta_0+2c'\zeta_0'+c\zeta_0'')-\frac{1}{r}(c'\zeta_0+c\zeta_0')=
\frac{1}{(1+\frac{r^2}{8})^2}c\zeta_0
\end{equation}
and because $\zeta_0$ is a solution of (\ref{7}) we get

\begin{equation}
\label{9}
-c''\zeta_0-c'(2\zeta_0'-\frac{1}{r}\zeta_0)=
0
\end{equation}
Setting $z=c'$ we obtain

\begin{equation}
\label{10}
z(r)=\frac{C}{r\zeta_0^2(r)}=C\frac{(8+r^2)^2}{r(8-r^2)^2}\sim\frac{C}{r}\quad\hbox{for $r$ large}
\end{equation}
where $C$ is a constant. This implies $c(r)\sim\log(r)$ for $r$ large. Hence $c\not\in L^\infty(\R)$
and $a\ fortiori$, $w\not\in L^\infty(\R)$. Then $\zeta_0(r)$ is the unique bounded solution of
(\ref{8}).

Now we consider the case $k=1$ in (\ref{6}). Here we have that
$\zeta_1(r)=\frac{r}{1+\frac{r^2}{8}}$ is a solution of (\ref{6}). Repeating the same argument as in
the case $k=0$ we obtain that a second linearly independent solution $w$ verifies

\begin{equation}
\label{11}
w(r)\sim r\quad\hbox{for $r$ large}.
\end{equation}
Hence again $w\not\in L^\infty(\R)$ and then $\zeta_1$ is the unique bounded solution of (\ref{6})
for $k\ge1$.

Now let us show that (\ref{6}) has no nontrivial solution for $k\ge2$.
For $k\ge1$ we set

\begin{equation}
\label{12}
A_k(\psi)=-\psi_k''-\frac{1}{r}\psi_k'+k^2\frac{\psi_k}{r^2}-
\frac{1}{(1+\frac{r^2}{8})^2}\psi_k.
\end{equation}
By contradiction let us suppose that there exists $\bar\psi\not\equiv0$ such that $A_k(\bar\psi)=0$
for some $k\ge2$. We claim that 

\begin{equation}
\label{13}
\bar\psi>0\quad\hbox{in }\R.
\end{equation}
Indeed if $\bar\psi$ changes sign we can select an interval $[x_1,x_2]$ with $0\le x_1<x_2<+\infty$
satisfying:

\begin{equation}
\label{14}
\bar\psi>0\quad\hbox{in }]x_1,x_2[.
\end{equation}
By (\ref{14}) we have that $\lambda_1(A_k)=0$ in $[x_1,x_2]$. On the other hand we have that
$A_k(\zeta_1)>0$ in $R$ and then the maximum principle holds in $[x_1,x_2]$ for $A_k$. Hence
$\lambda_1(A_k)>0$ and this gives a contradiction. Thus (\ref{13}) holds. Moreover we have that

\begin{equation}
\label{15}
\lim\limits_{r\rightarrow\infty}\bar\psi'(r)=0.
\end{equation}

In fact from (\ref{6})

\begin{equation}
\label{16}
(r\bar\psi')'-\frac{k^2}{r}\bar\psi=-\frac{r}{(1+\frac{r^2}{8})^2}\bar\psi
\end{equation}
and then, for $r>1$,

\begin{eqnarray}
\label{17}
& &
|r\bar\psi'(r)|\le|\bar\psi'(1)|+k^2\int_1^r\frac{|\bar\psi(t)|}{t}dt+\int_1^r\frac{t|\bar\psi(t)|}
{(1+\frac{t^2}{8})^2}dt\le\nonumber\\ 
& &\le \bar\psi'(1)+k^2||\bar\psi||_\infty\log
r+||\bar\psi||_\infty\int_1^\infty\frac{t} {(1+\frac{t^2}{8})^2}dt
\end{eqnarray}
and thus (\ref{17}) implies (\ref{15}).

Let us introduce the function $\eta(r)=r(\zeta_1\bar\psi'-\zeta_1'\bar\psi)$. It is easy to verify
that 
\begin{equation}
\label{18}
 \eta'(r)=(1-k^2)\frac{\bar\psi\zeta_1}{r}<0
\end{equation}
and
\begin{equation}
\label{19}
\lim\limits_{r\rightarrow\infty}\eta(r)=0.
\end{equation}
Thus, if we show that
$\lim\limits_{r\rightarrow0}\eta(r)=0$, using (\ref{18}) and (\ref{19}) we deduce a contradiction.
By the definition of $\eta$ we get, as $r\rightarrow0$,

\begin{equation}
\label{20}
\eta(r)=\frac{r^2\bar\psi'(r)}{1+\frac{r^2}{8}}-\bar\psi(r)\frac{r(1-\frac{r^2}{8})}
{(1+\frac{r^2}{8})^2}=\frac{r^2\bar\psi'(r)}{1+\frac{r^2}{8}}+o(r)
\end{equation}
and since (\ref{16}) implies $\lim\limits_{r\rightarrow0}r^2\bar\psi'(r)=0$ we have the claim.
So there exists no solution to (\ref{6}) as $k\ge2$.

Recalling that $Y_0(\theta)= constant$ and $Y_1(\theta)=x_1,x_2$, by (\ref{3}) we derive (\ref{2}).
\end{pf}
{\bf Proof of Theorem \ref{t:14}}
\\
By contradiction let us assume that there exist sequences $p_n\rightarrow\infty$
and $v_n\equiv v_{\l,p_n}\in H^1_0(\O)$, $v_n\not\equiv0$ satisfying 

\begin{equation}
\label{22}
\left\{
\begin{array}{cccc}
 -\Delta v_n+\l_n v_n =p_nu_{\l,p_n}^{p_n-1}v_n & \mbox{ in }\, \O \\ 
         v_n=0             & \mbox{ on }\,  \partial  \Omega 
\end{array}
\right.
\end{equation}

Since $\O$ satisfies the assumptions of the Gidas-Ni-Nirenberg Theorem we have that
$u_{\l,p_n}(0)=||u_{\l,p_n}||_\infty$ and 
$u_{\l,p_n}$ is even in $x_1$ and $x_2$. Thus we may assume that 

\begin{equation}
\label{23}
v_n(x_1,x_2)=v_n(-x_1,x_2)=v_n(x_1,-x_2).
\end{equation}
Set  $\tilde v_n(x)=v_n\Big(\frac{x}{\sqrt{p_n-1}||u_{\l,p_n}||_\infty^{\frac{p-1}{2}}}\Big)$ and
$u_n(x)=\frac{1}{||u_{\l,p_n}||_\infty}u_{\l,p_n}
\Big(\frac{x}{\sqrt{p_n-1}||u_{\l,p_n}||_\infty^{\frac{p-1}{2}}}\Big)$

We have that $\tilde v_n$ satisfies

\begin{equation}
\label{24}
\left\{
\begin{array}{cccc}
 -\Delta\tilde v_n+\frac{\l}{(p_n-1)||u_{\l,p_n}||_\infty^{p-1}}\tilde v_n
=\frac{p_n}{p_n-1}u_n^{p_n-1}\tilde v_n & \mbox{ in }\, \O_n \\ 
         v=0             & \mbox{ on }\,  \partial  \O_n
\end{array}
\right.
\end{equation}
with $\O_n=\sqrt{p_n-1}||u_{\l_n,p_n}||_\infty^{\frac{p-1}{2}}\,\O$.
Finally we set
\begin{equation}
\label{25}
z_n=\frac{\tilde v_n}{||\tilde v_n||_\infty}.
\end{equation}
Of course $z_n$ satisfies

\begin{equation}
\label{26}
\left\{
\begin{array}{cccc}
 -\Delta z_n+\frac{\l}{(p_n-1)||u_{\l,p_n}||_\infty^{p-1}}z_n
=\frac{p_n}{p_n-1}u_{\l,p_n}^{p_n-1}z_n & \mbox{ in }\, \O_n \\ 
|z_n|\le1& \mbox{ in }\, \O_n \\
         z_n=0             & \mbox{ on }\,  \partial  \O_n.
\end{array}
\right.
\end{equation}

We want to pass to the limit in  (\ref{26}). Since $u_{\l_n,p_n}$ is a
solution of $P_{\l,p_n}$ we get, computing $P_{\l,p_n}$ at $x=0$,
\begin{equation}
\label{27}
\l_n\le ||u_{\l,p_n}||_\infty^{p_n-1}.
\end{equation}
Then
$\frac{\l}{(p_n-1)||u_{\l,p_n}||_\infty^{p-1}}\le\frac{1}{p_n-1}\rightarrow0$ as
$n\rightarrow\infty$.
Now let us show the following estimates

\begin{equation}
\label{28}
\int_\O|\nabla z_n|^2 + \frac{\l}{(p_n-1)||u_{\l,p_n}||^{p_n-1}} \int_{\O_n}z_n^2 \le C_0,
\end{equation}
where $C_0$ is a positive constant independent of $n$.
Indeed from (\ref{26}) and (\ref{b0}) we derive 

\begin{equation}
\label{29}
\int_{\O_n}|\nabla z_n|^2 +  \frac{\l}{(p_n-1)||u_{\l,p_n}||^{p_n-1}} \int_{\O_n}z_n^2 \le 2\int_{\O_n}u_n^{p_n-1}=2(p_n-1)\int_{\O}u_{\l,p_n}^{p_n-1}
\le C_0.
\end{equation}
Moreover, from the classical Sobolev inequality
$$
\int_{\O_n} |\n \Phi|^2 + \l \int_\O |\Phi|^2 \geq C(\l,\O) \left(\int_\O |\Phi|^p\right)^{2/p}
$$
we deduce 
\begin{align*}
\int_{\O_n}|\nabla z_n|^2& + \frac{\l}{(p_n-1)||u_{\l,p_n}||^{p_n-1}} \int_{\O_n}z_n^2\\
& \geq C(\l,\O) \left(\frac{1}{(p_n-1)||u_{\l,p_n}||^{p_n-1}}\right)^{2/p_n} \left(\int_{O_n}|z_n|^{p_n}\right)^{2/p_n}
\end{align*}
From Lemma \ref{l:21} and \eqref{28} we get
$$
\left(\int_{O_n}|z_n|^{p_n}\right)^{2/p_n} \leq C,
$$
where $C$ is a constant independent of $n$.\\
Using (\ref{28}) and the standard regularity theory we deduce the existence of a function $z\in
C^2(\R^2)$, $|z|\le1$, such that $z_n\rightarrow z$ in $C^2_{loc}(\R^2)$. Moreover $z$ satisfies

\begin{equation}
\label{32}
\left\{
\begin{array}{cccc}
 -\Delta z=\frac{1}{(1+\frac{|x|^2}{8})^2}z & \mbox{ in }\R^2 \\
         |z|\le1& \mbox{ in }\R^2 \\
\int_{\R^2}|\nabla z|^2\le C_0
\end{array}
\right.
\end{equation}
From Lemma \ref{0} it follows that

\begin{equation}
\label{33}
z(x)=\sum_{i=1}^2a_i\frac{x_i}{1+\frac{|x|^2}{8}}+b\frac{8-|x|^2}{8+|x|^2}
\end{equation}
\\
{\bf  Step 1: $a_1=a_2=0$} in (\ref{33}). 

By (\ref{23}) we derive that $z(x)$ is even in $x_1$ and $x_2$. Hence by (\ref{33}) we deduce
$a_1=a_2=0$.
\\
{\bf  Step 2: $b=0$} in (\ref{33}).
From the previous step we have

\begin{equation}
\label{34}
z(x)=b\frac{8-|x|^2}{8+|x|^2}
\end{equation}

If $b\ne0$ we get that $\int_{\R^2}|\nabla z|^2=+\infty$ which is not possible. So $b=0$.
\\
{\bf  Step 3:} the contradiction.

In this step we prove the claim of Theorem \ref{t:14}. We point out that in this step we will use
Theorem \ref{t:13}.

By Step 1 and 2 we get that 

\begin{equation}
\label{34a}
z(x)\equiv0\quad\hbox{in }\R^2
\end{equation}

Since $||z_n||_{\infty}=1$ we can assume that there exists $x_n\in\O_n$ such that $z_n(x_n)=1$.
Since $z_n\rightarrow0$ in $C^2(\R^2)$ we obtain that $|x_n|\rightarrow\infty$.
By Theorem \ref{t:13} and Theorem \ref{t:12} we deduce that $u_{\l_n,p_n}^{p_n-1}(x_n)\rightarrow0$ as
$n\rightarrow\infty$. Otherwise, if by contradiction $u_{\l,p_n}^{p_n-1}(x_n) \geq C>0$ we derive the existence of a point $y_n$ such that $\n u_{\l,p_n}(y_n) = 0$, a contradiction with Theorem \ref{t:13}.\\
 Setting

\begin{equation}
\label{34b}
\bar z_n(x)=z_n(x+x_n)
\end{equation}
we get that $\bar z_n$ verifies

\begin{equation}
\label{35}
\left\{
\begin{array}{cccc}
 -\Delta \bar z_n+\frac{\l}{(p_n-1)||u_{\l,p_n}||_\infty^{p-1}}\bar z_n
=\frac{p_n}{p_n-1}u_n^{p_n-1}(x+x_n)\bar z_n & \mbox{ in }\, \O_n\setminus\{x_n\} \\ 
\bar z_n\le1& \mbox{ in }\, \O_n\setminus\{x_n\} \\
         \bar z_n(0)=1         \\
\int_{\R^2}|\nabla\bar z_n|^2\le C_0
\end{array}
\right.
\end{equation}

Passing to the limit in (\ref{35}) we derive that $\bar z_n\rightarrow\bar z$ in $C^2_{loc}(\R^2)$
where $\bar z$ satisfies

\begin{equation}
\label{36}
\left\{
\begin{array}{cccc}
 \Delta \bar z=0 & \mbox{ in } D \\ 
\bar|z|\le1& \mbox{ in } D,
\end{array}
\right.
\end{equation}
where $D$ is an half space if $dist(x_n,\O_n)\le K$ or $D=\R^2$ if
$\lim\limits_{n\rightarrow\infty}dist(x_n,\O_n)=+\infty$. In both case, by Liouville's Theorem we
have that
 
\begin{equation}
\label{36a}
z\equiv C\mbox{ in } D
\end{equation}
If $D$ is an half space, using that  $\bar z_n(x)=0$ for $x\in\partial\O_n-x_n$ we get
$\bar z(x)=0$
for $x\in\partial D$. Standard arguments (\cite{GS}) leads a contradiction with $\bar z_n(x_n)=1$.

Thus $D=\R^2$ and $\bar z(0)=1$. Moreover, since $\bar z_n\rightarrow1$ in $C^2(B(0,1))$, we get
$$
||z_n||_{L^{p_n}(\O_n)} \geq ||\bar z_n||_{L^{p_n}(B(0,1)} > \frac{1}{2}\quad\mbox{for }\,\, n>n_1.
$$
On the other hand since $z_n(x_n) =1$ and $z_n=0$ on the boundary of $\O_n$ we get that there exists a point $x_{2,n}$ with $|x_{2,n}|\to + \infty$ such that $z_n(x_{2,n})=\frac{1}{2}$. Seting
$$
\bar\bar{z}_n(x) = z_n(x+ x_{2,n})
$$
and repeating the same procedure of above we derive
$$
||z_n||_{L^{p_n}(\O_n)} \geq ||\bar\bar z_n||_{L^{p_n}(B(0,1))} > \frac{1}{2}\quad\mbox{for }\,\, n>n_2.
$$
Iterating this procedure, after a finite number of steps we reach a contradiction with \eqref{29}.\\ \\
{\bf Acknowledgments.} The main part of this work was done while the first author was
visiting the Department of Mathematics of the University of Roma " La
Sapienza" supported by the " Istituto Nazionale di Alta Matematica
". He would like to thank Mathematics Department for its warm
hospitality.

\end{document}